\documentclass[11pt,oneside]{amsart}
\usepackage{epsfig}
\usepackage{amsmath,ifthen, amsfonts, amssymb, srcltx, amsopn}
\usepackage{graphicx}
\usepackage[small]{caption}
\usepackage{enumerate}

\newcommand{\showcomments}{yes}
\renewcommand{\showcomments}{no}

\newsavebox{\commentbox}
\newenvironment{com}%
{\ifthenelse{\equal{\showcomments}{yes}}%
{\footnotemark
    \begin{lrbox}{\commentbox}
    \begin{minipage}[t]{1.25in}\raggedright\sffamily\tiny
    \footnotemark[\arabic{footnote}]}
{\begin{lrbox}{\commentbox}}}%
{\ifthenelse{\equal{\showcomments}{yes}}%
{\end{minipage}\end{lrbox}\marginpar{\usebox{\commentbox}}}
{\end{lrbox}}}

\newtheorem{thm}{Theorem}[section]
\newtheorem{lem}[thm]{Lemma}

\newtheorem{cor}[thm]{Corollary}

\newtheorem{prop}[thm]{Proposition}

\theoremstyle{definition}
\newtheorem{defn}[thm]{Definition}

\newtheorem{rem}[thm]{Remark}

\newcommand{\boundary}   {{\ensuremath \partial}}

\newtheorem*{Mainnonembedding}{Theorem 4.15}
\newtheorem*{Mainfreeint}{Theorem 4.11}

\begin{document}

\title[A Non-Quasiconvex Embedding of Relatively Hyperbolic Groups]
{A Non-Quasiconvex Embedding of Relatively Hyperbolic Groups}

\author[H.~Bigdely]{Hadi Bigdely}
      \address{Dept. of Math \& Stats.\\
               McGill University \\
               Montreal, Quebec, Canada H3A 2K6 }
\email{bigdely@math.mcgill.ca}


\date{\today}

\begin{com}
{\bf \normalsize COMMENTS\\}
ARE\\
SHOWING!\\
\end{com}

\begin{abstract}
For any finitely generated, non-elementary, torsion-free group $G$ that is hyperbolic relative to $\mathbb P$, we show that there exists a group $G^*$ containing $G$ such that  $G^*$ is hyperbolic relative to $\mathbb P$ and $G$ is not relatively quasiconvex in $G^*$. This generalizes a result of I. Kapovich for hyperbolic groups. We also prove that any torsion-free group $G$ that is non-elementary and hyperbolic relative to $\mathbb P$, contains a rank 2 free subgroup $F$ such that the group generated by ``randomly" chosen elements $r_1,\dots,r_m$ in $F$ is aparabolic, malnormal in $G$ and quasiconvex relative to $\mathbb P$ and therefore hyperbolically embedded relative to $\mathbb P$. 
\end{abstract}

\maketitle

\tiny
%\tableofcontents
\normalsize
\section{Introduction}
Gromov \cite{Gromov87} introduced relatively hyperbolic groups as a generalization of fundamental groups of complete finite-volume manifolds of pinched negative sectional curvature. Later Farb \cite{Farb98}, Bowditch \cite{Bowditch99RH} and Osin \cite{OsinBook06} formulated these groups with different viewpoints, however their notions are equivalent for finitely generated groups~\cite{HruskaRelQC}. In this paper, we follow Bowditch's approach. 

Relatively quasiconvex subgroups play an important role in the theory of relatively hyperbolic groups. For instance, any relatively quasiconvex subgroup is relatively hyperbolic and the intersection of two of these subgroups is again relatively quasiconvex. Moreover for a subgroup $G$ of Isom(X) where $X$ is a simply connected Riemannian manifold with pinched negative curvature, relatively quasiconvex subgroups of $G$ corresponds to ``geometrically finite" subgroups of $G$ (see \cite{HruskaRelQC}).

The main result of this paper is the following generalization of a result of Kapovich in \cite{Ikapovich99}. 

\begin{Mainnonembedding}\label{mainintro}
Let $G$ be a f.g., torsion-free group that is non-elementary and hyperbolic relative to $\mathbb P$. There exists a group $G^*$ that is hyperbolic relative to $\mathbb P$ such that $G$ is a subgroup of $G^*$ and $G$ is not quasiconvex in $G^*$ relative to $\mathbb P$.
\end{Mainnonembedding}

Let $G$ be hyperbolic relative to $\mathbb P$. We say that $G$ is \emph{absolutely relatively quasiconvex} if for any group $G^*$ that is hyperbolic relative to $\mathbb P$, and contains $G$ as a non-parabolic subgroup, $G$ is quasiconvex in $G^*$ relative to $\mathbb P$. Theorem 4.15 implies that any f.g., torsion-free and absolutely relatively quasiconvex group is infinite cyclic.

Let $G$ be hyperbolic relative to $\mathbb P$. As defined in \cite{Osin06}, a subgroup $H\leq G$ is \emph{hyperbolically embedded relative to $\mathbb P$}  if $G$ is hyperbolic relative to $\mathbb P\cup \{H\}$. We also prove the following result and note that the existence of hyperbolically embedded free subgroup in the following Theorem, has been shown in \cite{MatsudaOguniYamagata2012}.

\begin{Mainfreeint}
Let $G$ be torsion-free, non-elementary and hyperbolic relative to $\mathbb P$. Then there exists a rank 2 free subgroup $F$ of $G$ such that ``generically" any f.g. subgroup of $F$ is aparabolic, malnormal in $G$ and quasiconvex relative to $\mathbb P$, and therefore hyperbolically embedded relative to $\mathbb P$. In particular, $G$ contains a hyperbolically embedded rank 2 free subgroup.
\end{Mainfreeint}

Here by \emph{generically} we mean the following: Let $F_n$ be a rank $n$ free group. Let $\mathsf{N}(n,m,t)$ denote the number of $m$-tuples $(r_1,\dots,r_m)$ of cyclically reduced words in $F_n$ such that $|r_i|\leq t$ for each $i$. Moreover, let $\mathsf{N}_{\mathcal P}(n,m,t)$ be the number of such $m$-tuples such that $\langle r_1,\dots,r_m\rangle$ is  aparabolic, malnormal and quasiconvex in $G$ relative to $\mathbb P$, then $$\lim_{t\rightarrow\infty}\frac{\mathsf{N}_{\mathcal P}(n,m,t)}{\mathsf{N}(n,m,t)}=1.$$ 
In this paper all the ``generic" results are ``exponentially generic", i.e. the limit is exponentially approaches to $1$. 

This result is achieved through the following steps. Consider two hyperbolic elements $g$, $\bar g$ in $G$ such that $\langle g,\bar g\rangle$ is not cyclic. As $G$ is not elementary, such elements exist. We first prove that when $n$ is large enough, $F=\langle g^n,\bar g^n\rangle$ is a relatively quasiconvex, free subgroup. Finally with the aide of a variantion on a result of Arzhantseva, we prove that generically any finitely generated subgroup of $F$ is aparabolic, malnormal in $G$ and quasiconvex relative to $\mathbb P$.

In Section \ref{sec:quasigeod}, we give a criterion for a path in a $\delta$-hyperbolic space to be a quasigeodesic. Our argument extracts the hyperbolic case of a relative hyperbolic argument given in \cite{HsuWiseCubulatingMalnormal}. Section \ref{sec:Malnormal} contains a generic property for subgroups of free groups. In Section \ref{sec:main results}, we prove the main results.

\begin{com}add Acknowledgements\end{com}

\section{Quasigeodesics in hyperbolic spaces}\label{sec:quasigeod}
\begin{figure}[h]
\includegraphics[width=2.7 in]{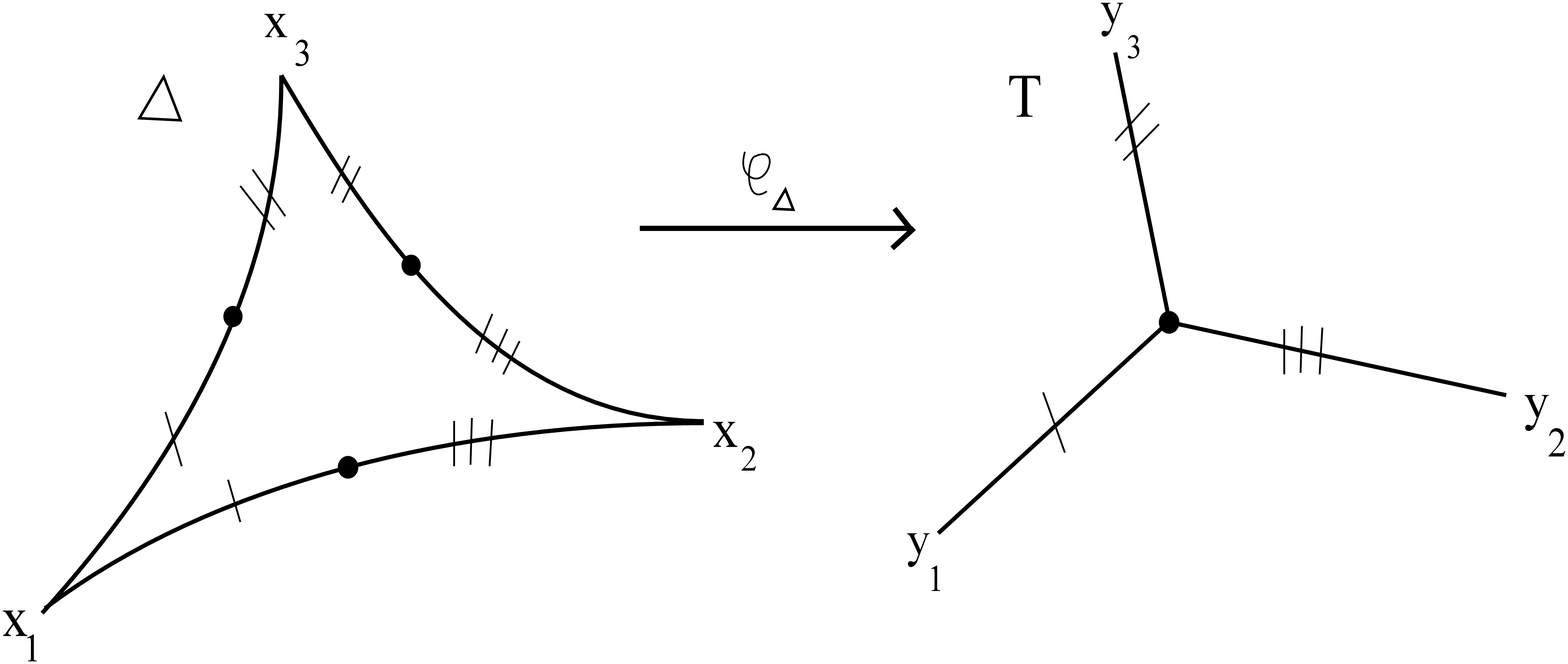}
\caption{}
\label{fig:triangle}
\end{figure}

%%%%%%%%%%%%%%%%%%%%%%%%%%%%%%%%%%%%%%%%%%%%%%%%%%%%%%%%%%%%%
A map $\varphi$ defined on a metric space $(X,d)$ is called \emph{$\epsilon$-thin} if $\varphi(x)=\varphi(y)$ implies $d(x,y)\leq \epsilon$. 
Let $X$ be a geodesic metric space and let $\bigtriangleup[x_1,x_2,x_3]$ be a geodesic triangle in $X$. Let $T$ be a tripod with three extremal vertices $y_1$, $y_2$ and $y_3$ so that $d(y_i,y_i)=d(x_i,x_j)$, see
Figure~\ref{fig:triangle}. The triangle $\bigtriangleup$ is called \emph{$\epsilon$-thin} if the map
$\varphi_{\bigtriangleup}: \bigtriangleup \rightarrow T$ which sends $x_i$ to $y_i$ and which is an isometry on the sides of $\bigtriangleup$, is
an $\epsilon$-thin map. A geodesic metric space $X$ is called \emph{$\delta$-hyperbolic}, for $\delta>0$, if any geodesic
triangle in $X$ is $\delta$-thin.

\begin{defn}[Quasigeodesic]
Let $[a,b]$ be a real interval and let $\lambda >0$ and $\epsilon \geq 0$. A $(\lambda,\epsilon)$-\emph{quasigeodesic} in a metric space $(X,d)$ is a function $\alpha : [a,b] \rightarrow X$ for all $s,t\in [a,b]$ such that $$d(s,t)\leq \lambda d(\alpha(s),\alpha(t))+\epsilon.$$  
\end{defn}
\begin{com}should I define fellow travel\end{com}
\begin{lem}\label{lem:Geodesic inequality}
Let $X$ be a $\delta$-hyperbolic space. Let $\tau_1$, $\sigma_2$, $\tau_2$, $\sigma_3$, $\tau_3$,$\dots$ be a sequence of concatenatable geodesics. Let $\gamma_k$ be a geodesic with the same endpoints of concatenation $\tau_1\sigma_2\tau_2\sigma_3\dots \sigma_k\tau_k$.
Suppose that there exists $c$ with the following properties:
\begin{enumerate}
\item \label{geod1}$|\tau_i|> 4c$ for $2\leq i \leq k-1$;
\item \label{geod2}$\tau_i$, $\tau_{i+1}$ can not $3\delta$-fellow travel for a distance $\geq c$;
\item \label{geod3}$\sigma_i$, $\tau_i$ and $\sigma_{i+1}$, $\tau_i$ can not $2\delta$-fellow travel for a distance $\geq c$.
\end{enumerate}
Then
\begin{enumerate}
\item \label{results of geod1} $|\gamma_{k+1}|\geq |\gamma_k|+|\sigma_{k+1}|+|\tau_{k+1}|-(6c+2\delta)$;
\item \label{results of geod2} The terminal subpaths of $\gamma_k$ and $\tau_k$, $\delta$-fellow travel for a distance of at least $|\tau_k|-2c$.
\end{enumerate}
\end{lem}

\begin{proof}
We prove the result by induction on $k$. When $k=1$, it is obvious. For each $k\geq 1$, let $\omega_k$ be a geodesic with the same endpoints as $\tau_1\sigma_2\tau_2\sigma_3\dots \sigma_k$. Consider two $\delta$-thin geodesic triangles $\gamma_k\sigma_{k+1}\omega_{k+1}$ and $\omega_{k+1}\tau_{k+1}\gamma_{k+1}$. Now consider the map to tripods with $\delta$-diameter fibers, and draw the geodesic triangles with sides $\leq \delta$. There are two cases according to the position of the points of these triangles on $\omega_{k+1}$. The first case is illustrated in Figure~\ref{fig:combined}-(II) and the second case in Figure~\ref{fig:quasigeodesic2}.

The notion $[i,j]$ is for the geodesic between points denoted by numbers $i$ and $j$ in the following figures and $|i,j|$ is the length of $[i,j]$.

First note that in both cases $4\in [2,5]$, otherwise if $4\in [1,2]$ by considering the point $13$ which is the comparison point of 2 on $\sigma_{k+1}$ in Figure~\ref{fig:combined}-(I), since the assumption~ \eqref{geod3} holds for $\sigma_{k+1}$, $\tau_k$, we have $|3,5|<c$. Also by induction the terminal subpaths of $\gamma_k$ and $\tau_k$, $\delta$-fellow travel for a distance of at least $|\tau_k|-2c$, therefore $|3,5|\geq |\tau_k|-2c$. Combining these two facts, $|\tau_k|<3c$ which is contradiction with assumption~\eqref{geod1}.

\begin{figure}[h]
\includegraphics[width=5.4 in]{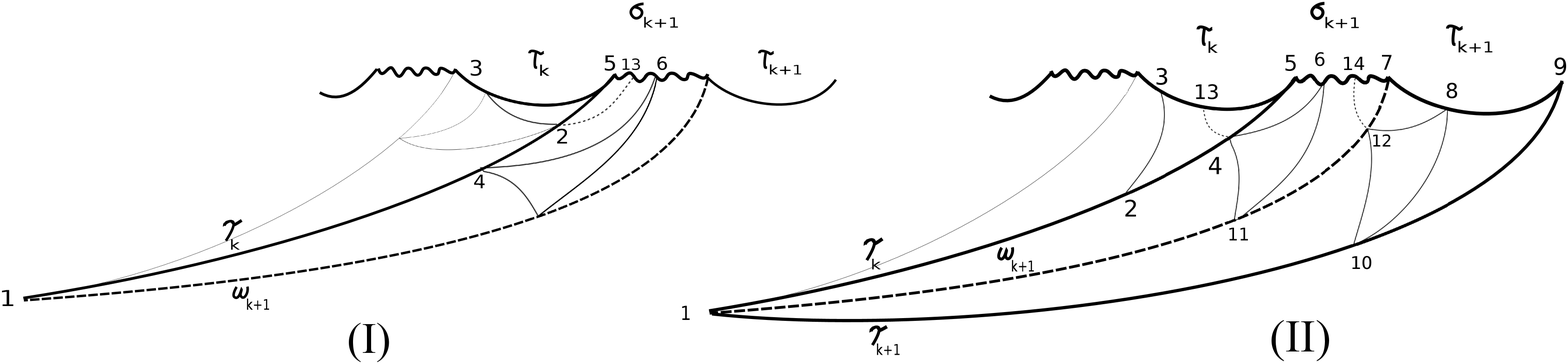}
\caption{}
\label{fig:combined}
\end{figure}

We consider the first case. Let the point $13$ be the comparison point of $4$ on $\tau_k$ and the point $14$ be the comparison point of $12$ on $\sigma_{k+1}$ and draw the geodesics $[4,13]$ and $[12,14]$.

%IN NIST\begin{figure}[h]
%\includegraphics[width=3 in]{quasigeodesic1.pdf}
%\caption{}
% \label{fig:quasigeodesic1}
%\end{figure}

Since by assumption~\eqref{geod3}, $\tau_k$ and $\sigma_{k+1}$ $2\delta$-fellow travel for a distance less than $c$, we have $|13,5|=|4,5|=|5,6|<c$. By similar argument $|14,7|=|7,12|=|7,8|<c$. Now
\begin{eqnarray}
|\omega_{k+1}| &=& |1,11|+|11,12|+|12,7|=|1,4|+|11,12|+|14,7|  \nonumber\\
& \geq & |1,4|+(|6,14|-2\delta)+|14,7|+(|4,5|-c)+(|5,6|-c) \nonumber \\
 & = & |1,5|+|5,7|-2\delta-2c=|\gamma_k|+|\sigma_{k+1}|-2\delta-2c  \nonumber
\end{eqnarray}

\begin{eqnarray}
|\gamma_{k+1}| &=& |1,10|+|10,9|=|1,12|+|8,9|=|1,11|+|11,12|+|8,9|  \nonumber\\
& \geq & |1,11|+|11,12|+|8,9|+(|12,7|-c)+(|7,8|-c) \nonumber \\
 & = & |1,7|+|7,9|-2c=|\omega_{k+1}|+|\tau_{k+1}|-2c  \nonumber\\
 &\geq &|\gamma_k|+|\sigma_{k+1}|+|\tau_{k+1}|-4c-2\delta \geq |\gamma_k|+|\sigma_{k+1}|+|\tau_{k+1}|-6c-2\delta  \nonumber
\end{eqnarray}

\begin{figure}[h]
\includegraphics[width=4 in]{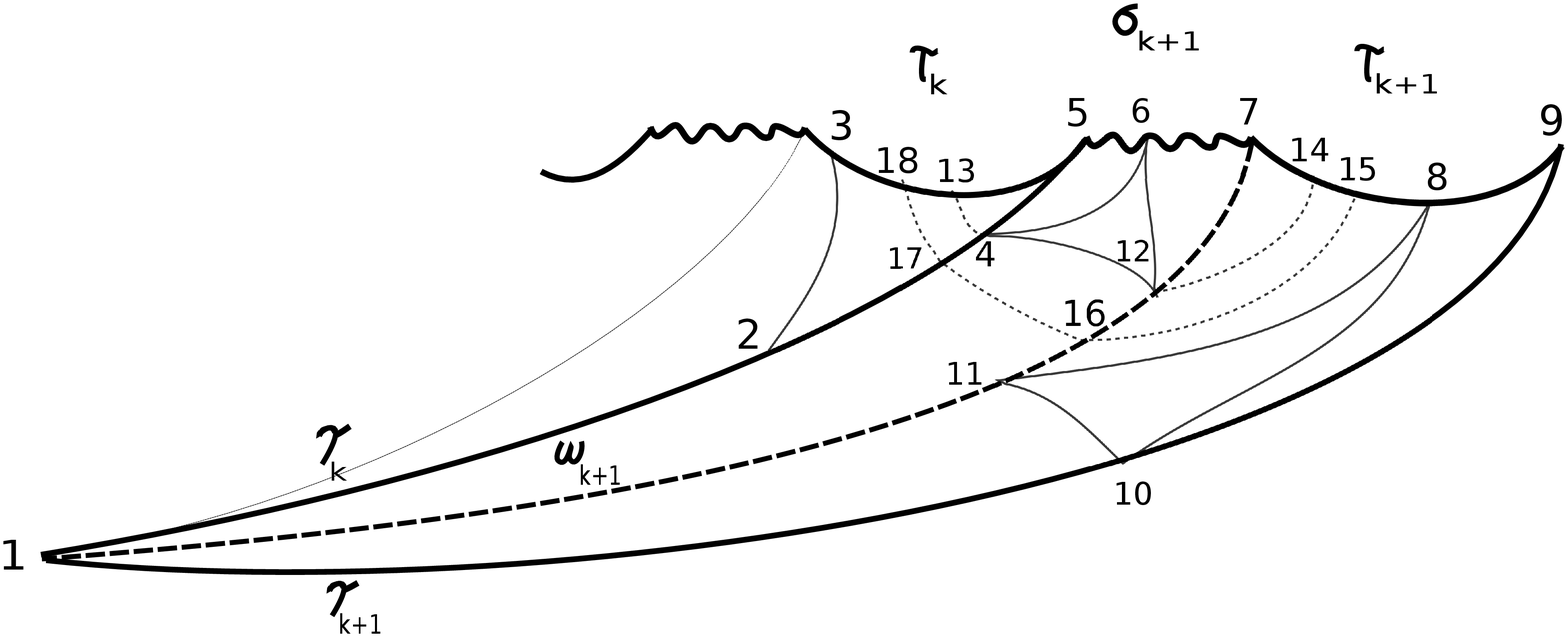}
\caption{}
\label{fig:quasigeodesic2}
\end{figure}

We now consider the second case. Note that the point $13$ is the comparison point of $4$ on $\tau_k$ and $14$ is the comparison point of $12$ on $\tau_{k+1}$. Assumption~ \eqref{geod3}, implies that $|13,5|=|4,5|=|5,6|<c$, by similar argument $|6,7|=|7,12|=|7,14|<c$.

First we prove $|7,8|$ is bounded by showing that $|7,8|<2c$. Assume that $|14,8|\geq c$, then we consider the point $15\in [14,8]$ such that $|14,15|=c$. Consider comparison points $16,17$ and $18$ of $15$ and geodesics between them. Note that $18\in [3,13]$, since by induction $|3,5|>2c$. Now since $|13,18|=|14,8|=c$ and by assumption~ \eqref{geod2}, $\tau_k$ and $\tau_{k+1}$ can not $3\delta$-fellow travel for a distance $\geq c$, we reach a contradiction. Therefore $|14,8|< c$ which implies $|7,8|<2c$.
Now we have:

\begin{eqnarray}
|\omega_{k+1}| &=& |1,12|+|12,7|=|1,4|+|6,7|  \nonumber\\
& \geq & |1,4|+|6,7|+(|4,5|-c)+(|5,6|-c) \nonumber \\
 & = & |1,5|+|5,7|-2c=|\gamma_k|+|\sigma_{k+1}|-2c  \nonumber
\end{eqnarray}

 \begin{eqnarray}
|\gamma_{k+1}| &=& |1,10|+|10,9|=|1,11|+|8,9|  \nonumber\\
& \geq & |1,11|+|8,9|+(|11,12|-c)+(|14,8|-c)+(|12,7|-c)+(|7,14|-c)\nonumber \\
 & = & |1,7|+|1,9|-4c=|\omega_{k+1}|+|\tau_{k+1}|-4c  \nonumber\\
 &\geq &|\gamma_k|+|\sigma_{k+1}|+|\tau_{k+1}|-6c \geq |\gamma_k|+|\sigma_{k+1}|+|\tau_{k+1}|-6c-2\delta  \nonumber
\end{eqnarray}\qedhere
\end{proof}

The following Lemma provides a criterion for a path in a $\delta$-hyperbolic space to be a quasigeodesic.
\begin{lem}[Quasigeodesic Criterion]\label{lem:quasigeodesic}
Let $X$ be a $\delta$-hyperbolic space. Let $$\gamma=\sigma_1\tau_1\sigma_2\tau_2\dots \sigma_n\tau_n\sigma_{n+1}$$ be a piecewise geodesic path. Suppose that:
\begin{enumerate}
\item \label{geod1}$|\tau_i|>2(6c+2\delta)$ for $i\geq 1$;
\item \label{geod2}$\tau_i$, $\tau_{i+1}$ can not $3\delta$-fellow travel for a distance $\geq c$;
\item \label{geod3}$\sigma_i$, $\tau_i$ and $\sigma_{i+1}$, $\tau_i$ can not $2\delta$-fellow travel for a distance $\geq c$.
\end{enumerate}
Then $\gamma$ is quasigeodesic. Indeed it is $(2,0)$-quasigeodesic.
\end{lem}

\begin{proof}
Let $\gamma'$ be a subpath of $\gamma$ and let $\lambda$ be a geodesic with the same endpoints as $\gamma'$. We show that $|\gamma'|\leq 2|\lambda|$. First note that any subpath $\gamma'$ of $\gamma$ can be expressed in the following form,
$$\gamma'= \overline{\tau_{i+1}}\sigma_{i+2}\tau_{i+2}\sigma_{i+3}\tau_{i+3}\dots \tau_{i+k-1}\sigma_{i+k}\overline{\tau_{i+k}}$$
where $\overline{\tau_{i+1}}$ and $\overline{\tau_{i+k}}$ are respectively subpathes of $\tau_{i+1}$ and $\tau_{i+k}$. Note that $\overline{\tau_{i+1}}$, $\overline{\tau_{i+k}}$, $\sigma_{i+2}$ and $\sigma_{i+k}$ can be trivial paths. Without loss of generality, we can shift the indices to make computation easier so that
$$\gamma'= \overline{\tau_{1}}\sigma_{2}\tau_{2}\sigma_{3}\tau_{3}\dots \tau_{k-1}\sigma_{k}\overline{\tau_{k}}$$
Let $\lambda$ be a geodesic with the same endpoints as $\gamma'$. Since $\gamma'$ satisfies in hypotheses of Lemma~\ref{lem:Geodesic inequality}, letting $\gamma_k$ in the lemma equal to $\lambda$ and using the hypothesis that $|\tau_i|>2(6c+2\delta)$ for $2\leq i\leq k-1$, we have:

 \begin{eqnarray}
|\lambda|=|\gamma_k| & \geq & \sum_{i=2}^{k-1}|\tau_i|+\sum_{i=2}^{k}|\sigma_i|+|\overline{\tau_{1}}|+|\overline{\tau_{k}}|-(k-1)(6c+2\delta)  \nonumber\\
& \geq &  \sum_{i=2}^{k-1}[|\tau_i|-(6c+2\delta)]+\sum_{i=2}^{k}|\sigma_i|+(6c+2\delta)+|\overline{\tau_{1}}|+|\overline{\tau_{k}}|  \nonumber \\
 &\geq & \frac{1}{2} \sum_{i=2}^{k-1}|\tau_i|+ \frac{1}{2} \sum_{i=2}^{k}|\sigma_i|   +\frac{1}{2} \left( |\overline{\tau_{1}}|+|\overline{\tau_{k}}| \right)=\frac{1}{2}|\gamma'|   \nonumber
\end{eqnarray}\qedhere
\end{proof}
%The following was proved in \cite[Thm~2T]{Tukia94}.\begin{com}check this\end{com}
%\begin{lem}\label{lem:characterizing subgps}
%Let $G$ be finitely generated and hyperbolic relative to $\mathbb P$. Let $H$ be non-elementary and non-parabolic subgroup of $G$ then $H$ contains an infinite hyperbolic element.
%\end{lem}
%The following was proved by Osin in \cite[Thm~4.3, Cor~1.7]{Osin06}.
%\begin{lem}\label{lem:maximal virtually cyclic}
%Let $G$ be hyperbolic relatively to $\mathbb P$. Let $g$ be a hyperbolic element of infinite order in $G$. Then $g$ is contained in a unique maximal elementary subgroup $E(g)$ of $G$. Also, $\langle g \rangle$ has finite index in $E(g)$ and $E(g)$ is hyperbolically embedded in $G$.
%\end{lem}
\section{A generic property for subgroups of free groups}\label{sec:Malnormal}
%%%%%%%%%%%%%%%%%%%%%%%%%%%%%%%%%%%%%%%%%%%%%%%%%%%%%%%%%%%%%
This section contains a variant result of Arzhantseva which says that given finitely many infinite index subgroups of a free group $F$, ``generically" subgroups of $F$ avoid conjugates of these infinite index subgroups. 
\begin{defn}[Almost Malnormal]
A subgroup $H$ is \emph{malnormal} in $G$ if $H\cap H^g=\{1\}$ for $g\notin H$,
 and similarly $H$ is \emph{almost malnormal} if this intersection $H\cap H^g$ is always finite.
 Likewise, a collection of subgroups $\{H_i\}$ is \emph{almost malnormal}
 if $H_i^g\cap H_j^h$ is finite unless $i=j$ and $gh^{-1}\in H_i$.
\end{defn}

\begin{defn}
Let $F_n$ be a nonabelian free group. Let $\mathsf{N}(n,m,t)$ denote the number of $m$-tuples $(r_1,\dots,r_m)$ of cyclically reduced words in $F_n$ such that $|r_i|\leq t$ for each $i$. Moreover, let $\mathsf{N}_{\mathcal P}(n,m,t)$ be the number of such $m$-tuples such that $\langle r_1,\dots,r_m\rangle$ has property $\mathcal P$. We say \emph{generically any subgroup of $F_n$ has property $\mathcal P$} if $$\lim_{t\rightarrow\infty}\frac{\mathsf{N}_{\mathcal P}(n,m,t)}{\mathsf{N}(n,m,t)}=1$$ 
\end{defn}

The following Proposition is a variant of Arzhantseva's result in \cite[Thm 1]{Arzhantseva2000}. We show the normalizer of a subgroup $H$ of $G$ by $\mathcal N(H)$.
\begin{prop}\label{prop:malnormal in virtually free}
Let $F_n$ be a nonabelian free group. Let $H_1,\dots,H_s$ be finitely generated infinite index subgroups of $F_n$. Generically, the group generated by randomly chosen words $r_1,\dots,r_m$ in $F_n$ has the property that $\langle r_1\dots,r_m\rangle\cap H_i^f=1$ for each $i$ and any $f\in F_n$. 
\end{prop}

\begin{proof}
Let $F_n=\langle x_1,\dots,x_n\rangle$. Let $\Gamma_{i}$ be the labelled directed graph corresponding to the core of the cover associated to $H_{i}$. $\Gamma_i$ is indeed the Stallings reduced folded graph. Let $(r_1,\dots,r_m)$ be an $m$-tuple of cyclically reduced words generated by $x_i^{\pm1}$, $1\leq i \leq n$ of length $|r_i|\leq t$.
By the argument in the proof of \cite[Thm 1]{Arzhantseva2000}, we have the following facts:

$\mathbf{(i)}$ The proportion of all $m$-tuples $(r_1,\dots,r_m)$ such that $\langle x_1,\dots,x_n~|~r_1,\dots,r_m\rangle$ is not $C'(\frac{1}{6})$ decreases exponentially when $t\rightarrow \infty$. So we assume that $\langle r_1,\dots,r_m\rangle$ satisfies the small cancellation condition $C'(\frac{1}{6})$.

$\mathbf{(ii)}$ We can assume that no cyclic shift of any $r_i$ contains a subword of length $\geq \frac{|r_i|}{2}$ which is the label of a path in some $\Gamma_{i}$, $1\leq i\leq s$. Since the proportion of $m$-tuples $(r_1,\dots,r_m)$ such that for at least one $k$, a subword of length $\geq \frac{|r_k|}{2}$ of $r_k$, is a label of a reduced path in some $\Gamma_{i}$, decreases exponentially when $t\rightarrow\infty$.

Now assuming (i)-(ii), we prove the result. It is enough to show that $\langle r_1,\dots,r_m\rangle^f \cap H_{i}=1$ for any $f\in F_n$ and each $i$. Assume that $\langle r_1,\dots,r_m\rangle^f\cap H_{i}\neq1$ for some $f\in F_n$ and some $i$. Then $\mathcal N(\langle r_1,\dots,r_m\rangle^f)\cap H_{i}\neq1$. Let $v\in \mathcal N(\langle r_1,\dots,r_m\rangle ^f)\cap H_{i}$. There is a nonempty reduced path in $\Gamma_{i}$ whose label is $v$. Also, $\bar v=v\mathcal N(\langle r_1,\dots,r_m\rangle ^f)$ is trivial element in the group $\frac{F_n}{\mathcal N(\langle r_1,\dots,r_m\rangle ^f)}$. Since $\frac{F_n}{\mathcal N(\langle r_1,\dots,r_m\rangle ^f)}$ satisfies $C'(\frac{1}{6})$, by Greendlinger's Lemma, (See~\cite{LS77}) there is some $r_j$ with a subword $v_j$ such that $|v_j|>\frac{|r_j|}{2}$ and, $r_j\cap v$ contains $v_j$. However, this contradicts (ii).
\end{proof}

Combining results in \cite[Thm~4.3]{BMNVW10} with Proposition~\ref{prop:malnormal in virtually free}, we get the following result:

\begin{cor}\label{cor:malnormal in virtually free}
Let $F_n$ be a nonabelian free group. Let $H_1,\dots,H_s$ be finitely generated infinite index subgroups of $F_n$. Exponentially generically, the group generated by randomly chosen words $r_1,\dots,r_m$ in $F_n$, is the free group of rank $m$, malnormal in $F_n$ and $\langle r_1\dots,r_m\rangle\cap H_i^f=1$ for each $i$ and any $f\in F_n$.
\end{cor}

\section{Main results}\label{sec:main results}
%%%%%%%%%%%%%%%%%%%%%%%%%%%%%%%%%%%%%%%%%%%%%%%%%%%%%%%%%%%%%
In this section we prove the main result. The following definition of relatively hyperbolic groups was given by Bowditch in \cite{Bowditch99RH}.

\begin{defn}[Relatively Hyperbolic] \label{defn:relatively hyperbolic}A {\em{circuit}}
in a graph is an embedded cycle. A graph $\Gamma$ is {\em{fine}} if each edge of $\Gamma$ lies in finitely many circuits of length $n$ for each $n$.

A group G is \emph{hyperbolic
relative to a finite collection of subgroups} $\mathbb{P}$ if $G$ acts cocompactly(without inversions)
on a connected, fine, hyperbolic graph $\Gamma$ with finite edge stabilizers, such that each element of $\mathbb P$ equals the stabilizer of a vertex of $\Gamma$, and moreover, each infinite vertex stabilizer is conjugate to a unique element of $\mathbb P$.
We refer to a connected, fine, hyperbolic graph $\Gamma$ equipped with
such an action as a \emph{$(G; \mathbb{P})$-graph}. Subgroups of G that are conjugate into
subgroups in $\mathbb{P}$ are \emph{parabolic}.
\end{defn}
The following definition of relatively quasiconvex subgroups of relatively hyperbolic groups is given in \cite{MartinezPedrozaWise2011}. It is also shown in \cite{MartinezPedrozaWise2011} that this definition is equivalent to Osin's definition for countable relatively hyperbolic groups.
\begin{defn}[Relatively Quasiconvex]\label{defn:rqc}
Let $G$ be hyperbolic relative to $\mathbb{P}$. A subgroup $H$
of $G$ is \emph{quasiconvex relative to $\mathbb{P}$} if for some (and hence any) $(G; \mathbb{P})$-graph $K$, there is a nonempty connected and quasi-isometrically embedded, $H$-cocompact subgraph $L$ of $K$. In the sequel, we sometimes refer to $L$ as a \emph{quasiconvex $H$-cocompact subgraph} of $K$.
\end{defn}

\begin{rem}\label{rem: everything for qc}
It is immediate from the Definition~\ref{defn:rqc} that in a relatively hyperbolic group, any parabolic subgroup is relatively quasiconvex, and any relatively quasiconvex subgroup is also relatively hyperbolic. In particular, the relatively quasiconvex subgroup $H$ is hyperbolic relative to the collection  $\mathbb P_H$ consisting of representatives of $H$-stabilizers of vertices of $L\subseteq K$. Note that a conjugate of a relatively quasiconvex subgroup is also relatively quasiconvex.

Let $G$ be hyperbolic relative to $\mathbb P_G$.
Suppose that $B$ is relatively quasiconvex in $G$,
and note that $B$ is then hyperbolic relative to $\mathbb P_B$. If $A\leq B$ is quasiconvex relative to $\mathbb P_B$ then $A$ is quasiconvex relative to $\mathbb P_G$.(see e.g. ~\cite[Lem~2.3]{BigdelyWiseAmalgams}).
\end{rem}

The following well-known property was proven in \cite{Bowditch99RH} (see \cite[Lem~2.2]{MartinezPedrozaWise2011}).
\begin{lem}[Almost Malnormal]\label{malnormal}
Let $G$ be hyperbolic relative to $\mathbb P$ then $\{P^g~|~P\in \mathbb P, g\in G \}$ is an almost malnormal collection of subgroups.
\end{lem}

Recall that, in a $\delta$-hyperbolic space $X$, for any hyperbolic element $g\in G$, we call any geodesic with two disjoint endpoints in $\boundary X$ which is $g$-invariant, \emph{axis} of $g$.

\begin{lem}[Parabolic or Hyperbolic]\label{lem:par or lox}
Let $G$ be hyperbolic relative to $\mathbb P$ and let $K$ be a $(G;\mathbb P)$-graph. Then for any infinite order $g\in G$, either $g$ is parabolic or $g^m$ has a quasiconvex axis in $K$ for some $m$.
\end{lem}
\begin{proof}
Let $d$ be the graph metric for $K$. By \cite[lem~2.2, 3.4]{Bowditch08}, either $g^m$ has a quasiconvex axis in $K$ for some $m$, or some (therefore any) $\langle g \rangle$-orbit is bounded. Now we show that if some $\langle g \rangle$-orbit is bounded then $g$ is parabolic. Let $x\in K^0$ and assume that $gx\neq x$. Let $\gamma$ be a geodesic between $x$ and $gx$ and let $e$ be an edge in $\gamma$. Now consider geodesic triangles $\bigtriangleup_n$ where for each $n$ the vertices of $\bigtriangleup_n$ are $x$, $gx$ and $g^nx$. Since the set $\{g^nx~|~n\geq 1\}$ is bounded and $K$ is fine and for each $n$ the triangle $\bigtriangleup_n$ contains $e$, the cardinality of the set $\{\bigtriangleup_n~|~n\geq 1\}$ is finite. Therefore there are $i> j$ such that $g^ix=g^jx$ which implies $g^{i-j}x=x$. Therefore $gx=x$ otherwise if stab(x)=$P$ then $g^{i-j}\in P\cap P^{g^{i-j}}$ but $g$ has infinite order which is contradiction since by Lemma~\ref{malnormal}, the intersection of two parabolic subgroups is finite. Lemma~\ref{malnormal}.
\end{proof}

Let $H$ be a subgroup of a group $G$. The commensurator subgroup of $H$  is \emph{Comm$(H)=\{g~|~[H:H\cap H^g]<\infty\}$}. The following Lemma was proved by Kapovich in \cite{Ikapovich99} in the special case when $G$ is hyperbolic.

\begin{lem}[Commensurator]\label{lem:comm}
Let $G$ be torsion-free and hyperbolic relative to $\mathbb P$. Let $F$ be a rank $\leq 2$ free subgroup of $G$ which is relatively quasiconvex and non-parabolic. Then Comm$(F)=F$.
\end{lem}
\begin{proof}
Since $F$ is relatively quasiconvex in $G$, by \cite[Thm 1.6]{HruskaWisePacking}, $F$ has finite index in Comm$(F)$. Moreover, since $G$ is torsion-free and $F$ is a free group, Comm$(F)$ is torsion-free that has a free subgroup of finite index. By the theorem of J.Stallings \cite{Stallings68}, this implies
that Comm$(F)$ is itself a free group. Now suppose that Comm$(F)\neq F$ , that is the index of $F$ in Comm$(F)$ is greater than $1$. By the Theorem II of ~\cite[ch 17]{Kurosh60}, this implies that the rank of Comm$(F)$ is strictly less than the rank of $F$.
However the rank of $F$ is $\leq 2$ and $F\leq \text{Comm}(F)$ which gives us a contradiction.
\end{proof}
\begin{defn}[Aparabolic]
Let $G$ be hyperbolic relative to $\mathbb P$. A subgroup $H$ is called \emph{aparabolic} if $H\cap P^g=1$ for any $P\in \mathbb P$ and $g\in G$.
\end{defn}
\begin{defn}[Quasi-isometry]
Let $(X,d_X)$ and $(Y,d_Y)$ be metric spaces. A map $\varphi:X\rightarrow Y$ is \emph{quasi-isometry} if for every $x_1$, $x_2\in X$,
$$\frac{1}{\kappa}d_X(x_1, x_2)-\epsilon \leq d_Y\big(\varphi(x_1),\varphi(x_2)\big)\leq \kappa d_X(x_1,x_2)+\epsilon,$$
for some constants $\kappa \geq 1$ and $\epsilon \geq 0$. We also call $\varphi$ a $(\kappa,\epsilon)$-quasi-isometry.
\end{defn}
The following was proved in \cite{MatsudaOguniYamagata2012}. Here, we give an independent and easy proof using Lemma~\ref{lem:quasigeodesic} and Bowditch's fine graph approach to relatively hyperbolic groups.

\begin{thm}\label{thm:free sub}
Let $G$ be torsion-free, non-elementary and hyperbolic relative to $\mathbb P$. Let $g$ and $\bar g$ be hyperbolic elements of $G$ such that $\langle g, \bar g\rangle$ is not cyclic. Then there exists $k$ such that for any $n\geq k$ the subgroup $F=\langle g^n,\bar g^n\rangle$ is free of rank 2, aparabolic and quasiconvex in $G$ relative to $\mathbb P$.
\end{thm}

\begin{proof}
Let $K$ be a $(G; \mathbb{P})$-graph. We use $[x,y]$ to denote a geodesic in $K$ from the vertex $x$ to the vertex $y$, and $|[x,y]|$ denotes the length of $[x,y]$. The open $\epsilon$-neighbourhood of $x$ is denoted by $N_{\epsilon}(x)$. %Let $R=R(K)$ be such that if $\alpha_1$ and $\alpha_2$ are two $(2,0)$-quasigeodesics with the same endpoints in $K$ then $\alpha_1\subseteq N_{\epsilon}(\alpha_2)$. 
Note that since $\langle g,\bar g \rangle$ is not cyclic, $\langle g^r,\bar g^s\rangle$ is not cyclic for any $r$, $s\geq 1$. Indeed if $\langle g^r,\bar g^s\rangle$ is cyclic for some $r$, $s$ then $g^i=\bar g^j$, for some $i$, $j\geq 1$, and so $g$ commutes with $\bar g^j$. This implies that $\langle \bar g^j\rangle \subseteq \langle \bar g\rangle \cap \langle \bar g \rangle^g$, hence $[\langle \bar g\rangle ~:~\langle \bar g\rangle \cap \langle \bar g \rangle^g]<\infty$. Therefore $g\in \text{Comm}(\langle \bar g\rangle)$ which equals $\langle \bar g\rangle$ by Lemma~\ref{lem:comm}. This contradicts that $\langle g,\bar g\rangle$ is not cyclic.

For sufficiently large $n$, we will construct a quasiconvex tree $L\subseteq K$ upon which $\langle g^{n},\bar g^{n}\rangle$ acts freely and cocompactly. By Lemma~\ref{lem:par or lox}, $h=g^{n_1}$ has a quasiconvex axis $\Upsilon$ for some $n_1$ and $\bar h=\bar g^{n_2}$ has a quasiconvex axis $\bar \Upsilon$ for some $n_2$. Choose $x_h$ and $x_{\bar h}$ to be vertices in $\Upsilon$ and $\bar \Upsilon$ such that $$d(x_h,x_{\bar h})=\min\{d(x,y)~|~x ~\text{is a vertex in} ~\Upsilon~ \text{and} ~y~ \text{is a vertex in}~ \bar \Upsilon\}.$$ Either $d(x_h,x_{\bar h})=0$ or $d(x_h,x_{\bar h})>0$. We give the proof in the second case and the first case is similar. Let $\sigma=[x_h,x_{\bar h}]$ be a geodesic. Choose $\mu>0$ such that both $\Upsilon\subseteq K$ and $\bar \Upsilon\subseteq K$ are $\mu$-quasiconvex, i.e. if a geodesic $\alpha$ in $K$ has endpoints on $\Upsilon$ then $\alpha \subseteq N_{\mu}(\Upsilon)$. Moreover, let $\Upsilon$ and $\bar \Upsilon$ be $(\kappa,\epsilon)$-quasiconvex where $\kappa\geq 1$.

Let $\tau=[x_h,h^mx_h]$ and $\bar \tau=[x_{\bar h},{\bar h}^mx_{\bar h}]$. Choose $m$ large enough that $\tau\nsubseteq N_{3\delta+4\mu+R}(\bar \tau)$ and $\tau\nsubseteq N_{3\delta+4\mu+R}(h\bar \tau)$. Note that such $m$ exists as otherwise $\Upsilon$ and $\bar \Upsilon$ are coarsely the same and so $\langle h,\bar h\rangle=\langle g^{n_1},\bar g^{n_2}\rangle$ would be cyclic which is impossible, as shown. Also as $\Upsilon$ and $\bar \Upsilon$ are quasiconvex we can choose $m$ large enough so that  $\tau\nsubseteq N_{2\delta+4\mu}(\sigma)$ and $\bar \tau\nsubseteq N_{2\delta+4\mu}(\sigma)$. Let $c=\max \{|\tau|,|\bar \tau|\}$ where $\tau$ and $\bar \tau$ satisfy all previously mentioned conditions. Without loss of generality, we can assume $m$ is large enough that $|\tau|>2(6c+2\delta)+1$ and $|\bar \tau|>2(6c+2\delta)+1$.

Let $F=\langle h^m,{\bar h}^m \rangle$ and let $L$ be the graph consisting of the union of all $F$-translates of the connected graph $\Upsilon \cup \bar \Upsilon \cup \sigma$ (see Figure~\ref{fig:tree}-(I)). Let $\gamma$ be the fundamental domain for the actin of $\langle h^m\rangle$ on $\Upsilon$ and $\bar \gamma$ be the fundamental domain for the actin of $\langle \bar h^m\rangle$ on $\bar \Upsilon$. Let $J$ be the graph obtained from $\gamma \sqcup \sigma \sqcup \bar \gamma$ by identifying the endpoints of $\gamma$ with the initial vertex of $\sigma$ and identifying the endpoints of $\bar \gamma$ with the terminal vertex of $\sigma$ (see Figure~\ref{fig:tree}-(II)). Note that $\pi_1 J\cong F$. Let $\tilde J$ be the universal cover of $J$ and consider $\pi_1 J$-equivariant map $\varphi: \tilde J\rightarrow L$ where the group $\pi_1 J$ maps to $F=\langle h^m,\bar h^m\rangle$. We show $\varphi$ is a quasi-isometry which implies that $L$ is a quasiconvex subgraph of $K$ and by Definition~\ref{defn:rqc}, this implies $F\leq G$ is relatively quasiconvex subgroup of $G$. We also prove that $\varphi$ is injective which shows that $L$ is a tree, therefore $F$ is free.

\begin{figure}[h]
\includegraphics[width=4.3 in]{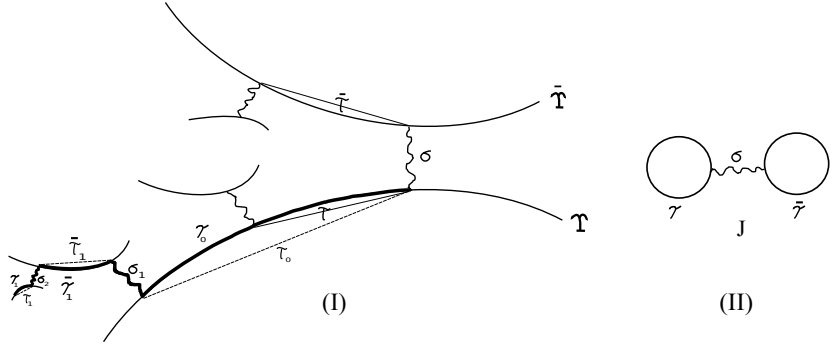}
\caption{The dark path is $\psi=\gamma_0\sigma_1\bar \gamma_1\sigma_2\gamma_1$. The path $\tau_i$ is a geodesic between endpoints of $\gamma_i$ and $\bar \tau_i$ is a geodesic between endpoints of $\bar\gamma_i$. 
 \label{fig:tree}}
\end{figure}

We show $\varphi$ is $\big(2(\kappa+\epsilon),0\big)$-quasi-isometry. Let $\xi$ be a path in $\tilde J$ and let $\varphi(\xi)=\psi$ such that $$\psi=\gamma_0\sigma_1\bar \gamma_1\sigma_2 \gamma_1\sigma_3\bar \gamma_2\sigma_4\gamma_2\cdots\sigma_r\gamma_r$$ where each $\sigma_i$ is an $F$-translate of $\sigma$, and each $\gamma_i$ and each $\bar\gamma_i$ is a subpath of an $F$-translate of $\Upsilon$ and $\bar \Upsilon$ respectively (see Figure~\ref{fig:tree}-(I)). Note that there are similar cases to consider where $\psi$ starts/ends with a proper path of $\sigma_i$ and where $\psi$ ends with $\bar \gamma_r$ instead of $\gamma_r$. 

Replace each $\gamma_i$ and each $\bar\gamma_i$ respectively by geodesics $\tau_i$ and $\bar\tau_i$ with the same endpoints to get $$\psi'=\tau_0\sigma_1\bar \tau_1\sigma_2\tau_1\sigma_3\bar \tau_2\sigma_4\tau_2\dots\sigma_r\tau_r$$
Now $\psi'$ satisfies the hypotheses of Lemma~\ref{lem:quasigeodesic}: Firstly, $|\tau_i|$, $|\bar \tau_i|>2(6c+2\delta)$ for each $i$ except perhaps for $i=0$ and $i=r$. Secondly, $\tau_i$ and $\bar \tau_i$ (and similarly $\bar \tau_i$ and $\tau_{i+1}$) do not $3\delta$-fellow travel for distance $\geq c$. Indeed $\tau_i\subseteq N_{\mu}(\gamma_i)$ and $\gamma_i\subseteq N_{\mu}(f\tau)$ where $f\tau$ is an $f$-translate of $\tau$ whose initial vertex is the same as $\gamma_i$. Therefore $\tau_i\subseteq N_{2\mu}(f\tau)$ and since $\tau\nsubseteq N_{3\delta+4\mu}(\bar \tau)$, we see that $\tau_i$ and $\bar\tau_i$ do not $3\delta$-fellow travel for distance $\geq c$. A similar argument shows that $\psi'$ satisfies the third property of Lemma~\ref{lem:quasigeodesic}. Therefore $\psi'$ is a $(2,0)$-quasigeodesic. We now show that $|\psi|\leq (\kappa+\epsilon)|\psi'|$ and since $\psi'$ is a $(2,0)$-quasigeodesic, this implies that $\varphi$ is $\big(2(\kappa+\epsilon),0\big)$-quasi-isometry. 

First note that since $\Upsilon$ and $\bar \Upsilon$ are $(\kappa,\epsilon)$-quasiconvex, $|\gamma_i|\leq \kappa |\tau_i|+\epsilon$ and $|\bar \gamma_i|\leq \kappa|\bar \tau_i|+\epsilon$, combining this with the assumption $|\tau_i|$ and $|\bar \tau_i|\geq 1$, we see that $|\gamma_i|\leq (\kappa+\epsilon) |\tau_i|$ and $|\bar \gamma_i|\leq (\kappa+\epsilon)|\bar \tau_i|$ (note that $\kappa+\epsilon\geq 1$). Therefore  

\begin{eqnarray}
|\psi| & = & |\gamma_0|+|\sigma_1|+|\bar \gamma_1|+|\sigma_2|+|\gamma_1|+|\sigma_3|+|\bar \gamma_2|+|\sigma_4|+|\gamma_2|+\dots+|\sigma_r|+|\gamma_r| \nonumber\\ & \leq & 
(\kappa+\epsilon)\big(|\tau_0|+|\sigma_1|+|\bar \tau_1|+|\sigma_2|+|\tau_1|+|\sigma_3|+|\bar \tau_2|+|\sigma_4|+|\tau_2|+\dots+|\sigma_r|+|\tau_r|\big) \nonumber\\
& = & (\kappa+\epsilon)|\psi'|\nonumber\
\end{eqnarray}
Thus $\varphi$ is quasi-isometry, therefore $L$ is a quasiconvex subgraph of $K$ and by Definition~\ref{defn:rqc}, this implies $F\leq G$ is relatively quasiconvex subgroup of $G$.

We now show that $\varphi$ is an injective map which implies $L$ is a tree and $F$ is free. Suppose $\varphi(x_1)=\varphi(x_2)$ where $x_1$ and $x_2$ are points on $\tilde J$. Since $\varphi$ is $\big(2(\kappa+\epsilon),0\big)$-quasi-isometry $$d_{\tilde J}(x_1,x_2)\leq 2(\kappa+\epsilon)d_K\big(\varphi(x_1),(x_2)\big)=0.$$ Therefore $x_1=x_2$ and $\varphi$ is injective.

Since $\varphi$ is an injective $F$-equivariant map and $F$ acts freely on $\tilde J$, $F$ acts freely on $L$. This implies that $F\cap P^g=1$ for any $g\in G$ and any $P\in \mathbb P$ and thus $F$ is aperiodic.
\end{proof}

Let $G$ be hyperbolic relative to $\mathbb P$. Recall that $Q\leq G$ is called \emph{hyperbolically embedded relative to $\mathbb P$} if $G$ is hyperbolic relative to $\mathbb P\cup \{Q\}$.

The following was proved in \cite{Osin06}.
\begin{thm}\label{hyp extend}
Let $G$ be hyperbolic relative to $\mathbb P$ and let $Q\leq G$ be a subgroup.Then $G$ is hyperbolic relative to $\mathbb P \cup \{Q\}$ if and only if $\mathbb P \cup \{Q\}$ is almost malnormal and $Q$ is quasiconvex in $G$ relative to $\mathbb P$
\end{thm}

\begin{lem}\label{thm:generically hyp}
Let $G$ be torsion-free, non-elementary and hyperbolic relative to $\mathbb P$. Let $F\leq G$ be a free and relatively quasiconvex subgroup that contains no non-trivial parabolic element. Then generically any f.g. subgroup of $F$ is aparabolic, malnormal in $G$ and quasiconvex relative to $\mathbb P$, and therefore hyperbolically embedded relative to $\mathbb P$.
\end{lem}
\begin{proof}
We show that generically any f.g. subgroup $M$ of $F$ is malnormal in $G$ and quasiconvex relative to $\mathbb P$. Clearly $M\cap P^g=1$ for any $g\in G$ and any $P\in \mathbb P$, thus $M$ is aparabolic. Therefore by Theorem~\ref{hyp extend}, $M$ will be hyperbolically embedded in $G$ relative to $\mathbb P$. Since $G$ is torsion-free and $F$ has rank 2, by Lemma~\ref{lem:comm}, Comm$(F)=F$. Since the intersection of $F$ with any parabolic subgroup is trivial, by~\cite[thm~1.6]{HruskaWisePacking}, $F$ has finite width i.e. there is a finite set of elements $\{g_1,\dots,g_s\}$ in $G-F$ such that if $F\cap F^g$ is infinite then $g\in g_i F$ for some $i$. Let $H_i=F\cap F^{g_i}$. Each $H_i$ is an infinite group which has infinite index in $F$. By Corollary ~\ref{cor:malnormal in virtually free}, generically, the subgroup $M=\langle r_1,\dots,r_m\rangle$ generated by randomly chosen words $r_1,\dots,r_m$ in $F$ is free group of rank $m$, malnormal in $F$ and $M\cap H_i^f=1$ for each $i$ and any $f\in F$. We show that $M$ is relatively quasiconvex and malnormal in $G$. Since $M$ is quasiconvex in $F$ and $F$ is relatively quasiconvex in $G$, by Remark ~\ref{rem: everything for qc}, $M$ is relatively quasiconvex in $G$.

We now show that $M$ is malnormal in $G$. Let $g\in G$ such that $M\cap M^g\neq 1$. If  $[F:F\cap F^g]<\infty$ then $g\in \text{Comm}(F)=F$, thus $g\in M$ since $M\leq F$ is malnormal. So we suppose that $[F:F\cap F^g]=\infty$, in which case $g=g_if$ for some $f\in F$ and some $i$. Now
$$1\neq M\cap M^g=M\cap M^{g_if}$$

therefore
$$1\neq M^{f^{-1}}\cap M^{g_i}~=~(M^{f^{-1}}\cap F)\cap M^{g_i}~=~M^{f^{-1}}\cap (M^{g_i}\cap F)~\subseteq ~M^{f^{-1}}\cap (F\cap F^{g_i})~=~M^{f^{-1}}\cap H_i$$
which contradicts that $M\cap H_i^f=1$ for any $f\in F$. Thus $M$ is malnormal in $G$
\end{proof}

Using Theorem ~\ref{thm:generically hyp} and Theorem~\ref{thm:free sub}, we get the following result. The existence of hyperbolically embedded free subgroup in the following Theorem, has been shown in \cite{MatsudaOguniYamagata2012}.
\begin{thm}\label{maincor:main}
Let $G$ be a torsion-free, non-elementary and hyperbolic relative to $\mathbb P$. Then there exists a rank 2 free subgroup $F$ of $G$ such that generically any f.g. subgroup of $F$ is aparabolic, malnormal in $G$ and quasiconvex relative to $\mathbb P$, and therefore hyperbolically embedded relative to $\mathbb P$. In particular, $G$ contains a hyperbolically embedded rank 2 free subgroup.
\end{thm}

\begin{proof}
 By \cite[Cor~4.5]{Osin06}, $G$ contains a hyperbolic element $g$. %By Lemma~\ref{lem:maximal virtually cyclic}, there is a unique maximal elementary subgroup $E(h)$ of $G$ where $\langle h \rangle$ has finite index in $E(h)$ and $E(h)$ is hyperbolically embedded in $G$.
By \ref{lem:comm} the commensurator subgroup of $\langle g \rangle$ is $\langle g\rangle$. Since $G$ is not elementary and not parabolic, there is a hyperbolic element $h$ in $G-\langle g \rangle$. Now the proof follows by Theorem ~\ref{thm:free sub} and Theorem ~\ref{thm:generically hyp}.
\end{proof}

The following Lemmas were proved in \cite{BigdelyWiseAmalgams}.
\begin{lem}[Combination]\label{thm:comb}
 Let $G_1$ be hyperbolic relative to $\mathbb P$. Let $P\in \mathbb P$ be isomorphic to ${P}'\leq P$. Let $G = G_1*_{{P}^t={P}'}$. Then $G$ is
hyperbolic relative to $\mathbb P-\{P\} \cup \{\langle P, t \rangle\}$.
\end{lem}

\begin{lem}[Quasiconvex Edges $\Longleftrightarrow$ Quasiconvex Vertices]\label{QC vertex}
Let $G$ be hyperbolic relative to $\mathbb P$. Suppose $G$ splits as a finite graph of groups whose vertex groups and edge groups are finitely generated. Then the edge groups are quasiconvex relative to $\mathbb P$ if and only if the vertex groups are quasiconvex relative to $\mathbb P$.
\end{lem}
The following was proved in \cite{DrutuSapir2005}.
\begin{lem}\label{associative hyp}
If $G$ is f.g. and hyperbolic relative to $\mathbb P=\{P_1,\dots, P_n\}$ and
each $P_i$ is hyperbolic relative to $\mathbb H_i=\{H_{i1},\dots,H_{i{m_i}}\}$, then $G$ is hyperbolic relative to $\bigcup_{1\leq i\leq n} \mathbb H_{i}$.
\end{lem}

The following Theorem generalizes Kapovich's result in \cite{Ikapovich99}.
\begin{thm}\label{main}
Let $G$ be a f.g. torsion-free group that is non-elementary and hyperbolic relative to $\mathbb P$. There exists a group $G^*$ that is hyperbolic relative to $\mathbb P$ such that $G$ is a subgroup of $G^*$ and $G$ is not quasiconvex in $G^*$ relative to $\mathbb P$.
\end{thm}

\begin{proof}
By Corollary~\ref{maincor:main}, there is a rank $2$ free subgroup $F=\langle a, b\rangle$ in $G$ such that $G$ is hyperbolic relative to $\mathbb P\cup \{F\}$. Let $\phi:F\rightarrow F$ such that $\phi(a)=abab^2\cdots ab^{100}$ and $\phi(b)=baba^2\cdots ba^{100}$. Then the group $$F^*=\langle a,b, t~|~ tat^{-1}=\phi(a), tbt^{-1}=\phi(b)\rangle.$$ is $C'(\frac{1}{6})$ and therefore hyperbolic. Moreover $F\leq F^*$ is exponentially distorted and thus $F$ is not quasiconvex in $F^*$. Indeed $d_{F}(1,t^nat^{-n})\geq 100^n$, where $d_F$ is the word metric for $F$. 

Consider the group $$G^*=\langle G,t~|~tat^{-1}=\phi(a), tbt^{-1}=\phi(b)\rangle.$$ We show that $G^*$ is hyperbolic relative to $\mathbb P$ and $G$ is not relatively quasiconvex in $G^*$. By Theorem~\ref{thm:comb}, since $G$ is hyperbolic relative to $\mathbb P\cup \{F\}$, the group $G^*$ is hyperbolic relative to $\mathbb P\cup\{\langle F,t\rangle\}$. Since $F^*=\langle F, t\rangle$ is hyperbolic, $G^*$ is hyperbolic relative to $\mathbb P$, by Lemma~\ref{associative hyp}.

Now suppose that $G$ is relatively quasiconvex in $G^*$. $$G^*=G*_{F^t=\varphi(F)}\cong G*_FF^*.$$
$F$ is relatively quasiconvex in $G$ and $G$ is relatively quasiconvex in $G^*$. Therefore $F$ is relatively quasiconvex in $G^*$. Therefore by Lemma~\ref{QC vertex}, $F^*$ is relatively quasiconvex in $G^*$. So $$F\leq F^*\leq G^*.$$ and both $F$ and $F^*$ are relatively quasiconvex in $G^*$, therefore $F$ is relatively quasiconvex in $F^*$. Since $F$ is aparabolic, and $F$ is relatively quasiconvex in $F^*$, we see that $F$ is quasiconvex in $F^*$ which is contradiction.
\end{proof}

\bibliographystyle{abbrv}
%\bibliographystyle{plain}

%\bibliography{C:/texfiles/LRQC/LRQC.tex/wise.bib}
%\bibliography{C:/papers/wise}
\bibliography{wise.bib}
\end{document}